\documentclass{amsart}

\usepackage{amsmath,amssymb,mathpazo}
\usepackage[all]{xy}
\SelectTips{cm}{10}

\newcommand {\g}{\mathfrak g}
\newcommand {\h}{\mathfrak h}
\renewcommand {\t}{\mathfrak t}
\newcommand {\n}{\mathfrak n}
\newcommand {\s}{\mathfrak s}
\renewcommand {\S} {\mathbf S}

\newcommand{\I}{\mathcal I}
\newcommand {\CC}{\mathbf C}
\newcommand {\F}{\mathbf F}

\newcommand {\CP}{\mathbf C P}

\newcommand{\Z}{\mathbf Z}

\newcommand{\op}{{\operatorname{op}}}

\newcommand{\id}{\operatorname{id}}

\newcommand {\sm}{\wedge}

\expandafter\def\csname c@figure\endcsname{\csname c@equation\endcsname}
\numberwithin{equation}{section}
\numberwithin{figure}{section}

\newtheorem{lemma}[equation]{Lemma}
\newtheorem{thm} [equation]{Theorem}
\newtheorem{corollary} [equation]{Corollary}
\newtheorem{prop} [equation]{Proposition}

\theoremstyle{definition}
\newtheorem{example}[equation]{Example}

\newtheorem*{defn}{Definition}

\newcommand{\hocolim}{\operatornamewithlimits{hocolim}}

\def\smashop#1_#2{%
\displaystyle{#1_{%
\hbox to 0pt{\hss$\scriptstyle{#2}$\hss}}\;}}
\DeclareMathOperator{\SU}{SU}

\makeatother
\newdir { >} {{}*!/-5pt/@{>}}

\DeclareMathOperator{\Top}{Top}
\newcommand{\spp}[1]{\S^0[#1]}
\DeclareMathOperator{\GL}{GL}

\newcommand{\bgtrans}{\overline\tau}

\author{Tilman Bauer}
\address{Institut f\"ur Mathematik\newline
\indent Westf\"alische Wilhelms-Universit\"at M\"unster\newline
\indent Einsteinstr. 62\newline
\indent 48149 M\"unster, Germany}
\email{tbauer@math.uni-muenster.de}
\author{Nat\`alia Castellana}
\address{Departament de Matem\`atiques\newline
\indent Universitat Aut\`onoma de Barcelona\newline
\indent 08193 Bellaterra, Spain}
\email{natalia@mat.uab.es}

\date{February 20, 2006}
\title{Adjoint spaces and flag varieties of $p$-compact groups}
\thanks{N. Castellana is partially supported by MCYT grant MTM2004-06686}
\thanks{Support by the Institut Mittag-Leffler (Djursholm, Sweden) is gratefully
acknowledged.}

\begin{document}

\begin{abstract}
For a compact Lie group $G$ with maximal torus $T$, Pittie and Smith showed that the flag variety $G/T$ is always a stably framed boundary. We generalize this to the category of $p$-compact groups, where the geometric argument is replaced by a homotopy 
theoretic argument showing that the class in the stable homotopy groups of spheres represented by $G/T$ is trivial, even $G$-equivariantly. As an application, we consider an unstable construction of a $G$-space mimicking the adjoint representation 
sphere of $G$ inspired by work of the second author and Kitchloo. This construction stably and $G$-equivariantly  splits off its top cell, which is then shown to be a dualizing spectrum for $G$.
\end{abstract}

\maketitle

\section{Introduction}

Let $G$ be a compact, connected Lie group of dimension $d$ and rank $r$ with maximal torus $T$. Left translation by elements of $G$ on the tangent space $\g = T_eG$ induces a framing of $G$. By the Pontryagin-Thom construction, $G$ with this framing 
represents an element $[G]$ in the stable homotopy groups of spheres; this has been extensively studied for example in \cite{Smith,Wood,Knapp,Ossa}.

The following classical argument shows that the flag variety $G/T$, while not necessarily framed, is still a stably framed manifold: since every element in a compact Lie group is conjugate to an element in the maximal torus, the conjugation map $G 
\times T \to G$, $(g,t) \mapsto gtg^{-1}$ is surjective, and furthermore, it factors through $c\colon G/T \times T \to G$. An element $s \in T$ is called \emph{regular} if the centralizer $C_G(s) \supseteq T$ equals $T$, or, equivalently, if 
$c|_{G/T \times \{s\}}$ is an embedding; it is a fact from Lie theory that the set of irregular elements has positive codimension in $T$. Thus there is a regular element $s$ such that the derivative of $c$ has full rank along $G/T \times \{s\}$, 
and by the tubular neighborhood theorem, it induces an embedding of $G/T \times U$, where $U$ is a contractible neighborhood of $s$ in $T$. Thus the framing of $G$ can be pulled back to a stable framing of $G/T$.

Pittie and Smith showed in \cite{Pittie:flag,Pittie-Smith:flag} that $G/T$ is always the boundary of another framed manifold $M$, and moreover, that $M$ has a $G$-action which agrees with the standard $G$-action on $G/T$ on the boundary. In terms 
of homotopy theory, this is saying that the class $[G/T] \in \pi_{d-r}^s$ induced by the Pontryagin-Thom construction is trivial.

The first main result of this paper generalizes this fact to $p$-compact groups.

\begin{thm} \label{pittiethm}
Let $G$ be a $\Z/p$-local, $p$-finite group with maximal torus $T$ such that $dim(G)>dim(T)$. Then the Pontryagin-Thom construction $[G/T]\colon S_G \to S_T$ is $G$-equivariantly null-homotopic, with $G$ acting trivially on $S_T$.
\end{thm}

The statement of this theorem requires some explanation. A \emph{$p$-compact group} \cite{dw} is a triple $(G,BG,e)$ such that
\begin{itemize}
\item $G$ is \emph{$p$-finite}, i.~e., $H_*(G;\F_p)$ is finite;
\item $BG$ is \emph{$\Z/p$-local}, i.~e. whenever $f\colon X \to Y$ is a mod-$p$ homology equivalence of CW-complexes, then $f^*\colon [Y,BG] \to [X,BG]$ is an isomorphism;
\item $e\colon G \to \Omega BG$ is a homotopy equivalence.
\end{itemize}

Clearly, $G$ and $e$ are determined by $BG$ up to homotopy, making this definition somewhat redundant. Although a priori $G$ is only a loop space, we will henceforth assume we have chosen a rigidification such that $G$ is actually a topological 
group. This is always possible, for example by using the geometric realization of Kan's group model of the loops on a simplicial set \cite{Kan}.

A $\Z/p$-local, $p$-finite loop space is only slightly more general than a $p$-compact group in that the latter also requires $\pi_0(G)$ to be a $p$-group. We will have no need to assume this in Theorem \ref{pittiethm}.

By \cite{dw}, every $p$-compact group has a maximal torus $T$; that is, there is a monomorphism $T \to G$ with $T \simeq L_p (\S^1)^r$ and $r$ is maximal with this property. By definition, a \emph{monomorphism} of $p$-compact groups is a group 
monomorphism $H \to G$ such that $G/H$ is $p$-finite (see \cite{bauer:pcfm} for this slightly nonstandard point of view). Dwyer and Wilkerson show that $T$ is essentially unique. Since a maximal torus is always contained in the 
identity component of a $p$-compact group, the same works for $\Z/p$-local, $p$-finite groups.

Denote by $\spp{X}$ the suspension spectrum of a space $X$ with a disjoint base point added. 

\begin{defn}[\cite{kl00}]
Let $G$ be a topological group. Define $S_G$, the \emph{dualizing spectrum of $G$}, to be the spectrum of homotopy fixed points of the right action of $G$ on its own suspension spectrum. That is, $S_G = (\spp{G})^{hG^\op}$ as left $G$-spectra.
\end{defn}

In \cite{bauer:pcfm}, the first author showed that for a connected, $d$-dimensional $p$-com\-pact group $G$, $S_G$ is always homotopy equivalent to a $\Z/p$-local sphere of dimension $d$. Furthermore, there is a $G$-equivariant logarithm map 
$\spp{G} \to S_G$, where $G$ acts on the left by conjugation. If $G$ is the $\Z/p$-localization of a connected Lie group, then $S_G$ is canonically identified with the suspension spectrum of the one-point compactification of the Lie algebra of $G$. 
Thus we may call $S_G$ the \emph{adjoint (stable) sphere} of $G$.

In a spectacular case of shortsightedness,  \cite{bauer:pcfm} restricts its scope to connected $p$-compact groups where everything would have worked for $\Z/p$-local, $p$-finite groups $G$ as well. In this case, $S_G$ has the mod-$p$ homology of a 
$d$-dimensional sphere. Similarly, the proof of the following was given in \cite[Cor. 24]{bauer:pcfm} for connected groups, but immediately generalizes.

Let $DM$ be the Spanier-Whitehead dual of a finite CW-spectrum $M$.
\begin{lemma} \label{relativeduality}
Let $H<G$ be a monomorphism of $\Z/p$-local, $p$-finite groups. Then there is a relative $G$-equivariant duality weak equivalence
\[
G_+ \sm_H S_H \simeq D\left(\spp{G/H}\right) \sm S_G.
\]
\end{lemma}

For any space $X$, there is a canonical map $\epsilon\colon \spp{X} \to \S^0$ given by applying the functor $\spp{-}$ to $X \to *$. If $T < G$ is a sub-torus in a $\Z/p$-local, $p$-finite group then there is a stable $G$-equivariant map
\begin{multline} \label{tpconstr}
[G/T]\colon S_G \xrightarrow{\id \sm D\epsilon} S_G \sm D\left(\spp{G/T}\right)\\ \underset{\text{Lemma }\ref{relativeduality}}\simeq G_+ \sm_T S_T \simeq \spp{G/T} \sm S_T \xrightarrow{\epsilon} S_T
\end{multline}
where the homotopy equivalence on the right hand side holds because $S_T$ has a homotopy trivial $T$-action as $T$ is homotopy abelian. The first map is studied in \cite{bauer:pcfm}. This is the map referred to in Theorem \ref{pittiethm}; it 
generalizes the Pontryagin-Thom construction.

\bigskip

In the second part of this paper, as an application of Theorem \ref{pittiethm}, we study the relationship between two notions of adjoint objects of $p$-compact groups. It is an interesting question to ask whether the action of $G$ on $S_G$ actually 
comes from an unstable action of $G$ on $\S^d$. We will not be able to answer this question in this paper. However, there is an alternative, unstable construction of an adjoint object for a connected $p$-compact group $G$ inspired by the following:

\begin{thm}[\cite{castellana-kitchloo:adjoint,mitchell:buildings}] \label{ckthm}
Let $G$ be a semisimple, connected Lie group of rank $r$. There exist subgroups $G_I < G$ for every $I \subsetneq \{1,\dots,r\}$ and a homeomorphism of $G$-spaces
\[
A_G := \Sigma \hocolim_{I \subsetneq \{1,\dots,r\}} G/G_I \to \g \cup \{\infty\}
\]
to the one-point compactification of the Lie algebra $\g$ of $G$.
\end{thm}

In the second part of this paper, we define a $G$-space $A_G$ for every connected $p$-compact group $G$ and show:

\begin{thm} \label{adjointcomparison}
For any connected $p$-compact group $G$, there is a $G$-equivariant splitting $\spp{A_G} \simeq S_G \vee R$ for some finite $G$-spectrum $R$.
\end{thm}

This result links the two notions of adjoint objects together. Thus stably, the adjoint sphere is a wedge summand of the adjoint space. 

Unfortunately, $A_G$ is in general not a sphere.

\subsection*{Acknowledgements} We would like to thank the Institut Mittag-Leffler for its support while finishing this work, and Nitu Kitchloo for helpful discussions.

\section{The stable $p$-complete splitting of complex projective space}

\subsection{Stable splittings from homotopy idempotents}

Let $p$ be a prime. We denote by $L_p$ the localization functor on topological spaces with respect to mod-$p$ homology, which coincides with $p$-completion on nilpotent spaces \cite{yellowmonster}. Let $S = L_p \S^1$ be the $p$-complete $1$-sphere, 
and set $P = \spp{BS}$.
It is a classical result that
\begin{equation} \label{splittingbs1}
P \simeq \bigvee_{s=0}^{p-2} P_s
\end{equation}
for certain $(2i-1)$-connected spectra $P_i$. In this section, we will investigate this splitting and its compatibility with certain transfer maps.

Let $X$ be a spectrum, $e \in [X,X]$ and define
\[
eX = \hocolim\{X \xrightarrow{e} X \xrightarrow{e} \cdots \}.
\]
If $e$ is idempotent, this is a homotopy theoretic analog of the image of $e$. Any such idempotent $e$ yields a stable splitting $X \simeq eX \vee (1-e) X$. If $\{e_1,\dots,e_n\}$ are a complete set of orthogonal idempotents (this means that each 
$e_i$ is idempotent, $e_ie_j \simeq *$, and $\id_X \simeq e_1+\cdots+e_n$), then they induce a splitting $X \simeq e_1 X \vee \cdots \vee e_n X$.

\begin{example}
Let $p$ be an odd prime. Denote by $\psi \colon P \to P$ the map induced by multiplication with a $(p-1)$st root of unity $\zeta$. Define $e_s\colon P \to P$ by
\[
e_s = \frac1{p-1}\left(\sum_{i=0}^{l-1} \zeta^{-is} \psi^i\right).
\]
It is straightforward to check that $\{e_0,\dots,e_{p-2}\}$ are a complete set of orthogonal idempotents in $[P,P]$. They induce the splitting \eqref{splittingbs1} by defining $P_s = e_s P$.

Setting $H_*(P) = \Z_p\{x_j\}$ with $|x_j|=2j$, we have that $(e_i)_*\colon H_*(P) \to H_*(P)$ is given by
\begin{equation} \label{splittingonhomology}
(e_i)_*(x_j) = \begin{cases} x_j; & j \equiv i \pmod{p-1}\\ 0; & \text{otherwise.} \end{cases}
\end{equation}
\end{example}

\subsection{Transfers as splittings} \label{transfersplitting}

Let $1 \to H \xrightarrow{i} G \to W \to 1$ be an extension of compact Lie groups. Then associated to the fibration $W \to BH \to BG$ there are two versions of functorial stable transfer maps \cite{Becker-Gottlieb:transfer-bundles, 
Becker-Gottlieb:transfer-fibrations}:
\begin{enumerate}
\item The Becker-Gottlieb transfer $\bgtrans\colon \spp{BG} \to \spp{BH}$ 
\item The stable Umkehr map $\tau\colon BG^\g \to BH^\h$ of Thom spaces of the adjoint representation of the Lie groups.
\end{enumerate}
Both versions can be generalized to a setting where the groups involved are not Lie groups but only $\Z/p$-local and $p$-finite \cite{dwyer,bauer:pcfm}. For such a group $G$, $BG^\g$ is defined to be the homotopy orbit spectrum of $G$ acting on the 
dualizing spectrum $S_G$; since $H_*(S_G) = H_*(\S^d;\Z_p)$, we have a (possibly twisted) Thom isomorphism between $H_*(BG)$ and $H_*(BG^\g)$.

Note that $\bgtrans$ factors through $\tau$ in the following way:
\begin{multline} \label{transferfactor}
\spp{BG} \xrightarrow{\tau'} BH^\nu \xrightarrow{\operatorname{comult.}} BH^\nu \sm_{BG} \spp{BH}\\ \xrightarrow{\id \sm \Delta} BH^\nu \sm_{BG} \spp{BH} \sm_{BG} \spp{BH} \xrightarrow{\operatorname{eval} \sm \id} \spp{BH}
\end{multline}
where $\nu = \h - i^*\g$ is the normal fibration along the fibers of $BH \to BG$, $\tau'$ is $\tau$ twisted by $-\g$, and the right hand side evaluation map is defined by identifying $BH^\nu$ with the fiberwise Spanier-Whitehead dual of $BH$ over 
$BG$. 

\begin{prop} \label{bgfactorization}
Let $W=C_l$ be a finite cyclic group acting freely on $S$, with $l \mid p-1$. Denote by $N = S \rtimes W$ the semidirect product with respect to this action. Then the Becker-Gottlieb transfer map $\bgtrans\colon \spp{BN} \to P$ factors through $f P 
\to P$ for some idempotent $f\colon P \to P$ which induces the same map in homology as $e_0 + e_l + \cdots + e_{p-1-l}$, and the induced map $\spp{BN} \to f P$ is a mod-$p$ homology equivalence.
\end{prop}

\begin{proof}
Since $p \nmid |W|$, the Serre spectral sequence associated to the group extension $S \xrightarrow{i} N \to W$ is concentrated on the vertical axis and shows that 
\[
H^*(BN;\Z_p) \cong H^*(BS;\Z_p)^W \cong \Z_p[z^l] \hookrightarrow \Z_p[z] \cong H^*(BS;\Z_p).
\]
In this case, the Becker-Gottlieb transfer is nothing but the usual transfer for finite coverings, therefore $i \circ \bgtrans$ is multiplication by $|W| = l \in \Z_p^\times$. Setting $I = l^{-1} i\colon P \to \spp{L_p BN}$, we thus get orthogonal 
idempotents in $[P,P]$:
\[
f = \bgtrans \circ I \quad \text{and} \quad e = \id_P - f.
\]

Clearly, $e \circ \bgtrans \simeq *$, thus $\bgtrans$ factors through $fP$ and induces an isomorphism $\spp{L_p BN} \to fP$, in particular a mod-$p$ homology isomorphism between $\spp{BN}$ and $fP$. The computation of the homology of $BN$ together 
with \eqref{splittingonhomology} implies that $f_* = (e_0 + e_l + \cdots + e_{p-1-l})_*$.
\end{proof}

\begin{corollary} \label{transferfactorization}
Let $S,\; N, \; W$ be as above. Then the stable Umkehr map
\[
BN^\n \to BS^\s \simeq \Sigma P
\]
factors through $\Sigma f P \to \Sigma P$ for some $f\colon P \to P$ which induces the same morphism in homology as $\sum_{i=0}^{\frac {p-1}l} e_{(i+1)l-1}$.
The induced map $BN^\n \to \Sigma f P$ is a mod-$p$ homology equivalence.
\end{corollary}
\begin{proof}
This follows from a similarly simple homological consideration.
The $S$-fibration $\n$ is not orientable, thus we have a twisted Thom isomorphism
\[
\tilde H^{n+1}(BN^\n) \cong H^n(BN;\mathcal H^1(S;\Z_p))
\]
where $\pi_1(BN) = \Z/l$ acts on $H^1(S;\Z_p) \cong \Z_p$ by multiplication by an $l$th root of unity.
Thus
\[
H^i(BN^\n;\Z_p) = \begin{cases} \Z_p; & i \equiv -1 \pmod{l}\\
0; & \text{otherwise.} \end{cases}
\]

The factorization \eqref{transferfactor} of  $\bgtrans$ through $\tau$ 
\begin{multline*}
\spp{BN} \xrightarrow{\tau'} BS^\nu \xrightarrow{\operatorname{comult.}} BS^\nu \sm_{BN} \spp{BS}\\ \xrightarrow{\id \sm \Delta} BS^\nu \sm_{BN} \spp{BS} \sm_{BN} \spp{BS} \xrightarrow{\operatorname{eval} \sm \id} \spp{BS}
\end{multline*}
simplifies considerably since $i^*\nu$ is the trivial $1$-dimensional fibration over $BS$, and the composition of the three right hand side maps is an equivalence.

In Prop.~\ref{bgfactorization} it was shown that $I \circ \bgtrans = \id_{\spp{L_pBN}}$, thus the same holds after twisting with $\n$:
\[
\id_{BN^\n}\colon L_p BN^\n \xrightarrow{\tau} BS^\s \to BS^{i^*\n} \xrightarrow{I^\n} L_p BN^\n.
\]
If we denote the composition $BS^\s \to BS^{i^*\n} \xrightarrow{I^\n} L_p BN^\n$ by $I$, overriding its previous meaning, the argument now proceeds as in Prop. \ref{bgfactorization}. Using the computation of $H^*(BN^\n;\Z_p)$, we find that $L_p 
BN^\n \simeq (\tau \circ I)P$, and 
\[
(\tau \circ I)_*=\sum_{i=0}^{\frac {p-1}l} (e_{(i+1)l-1})_*
\]
\end{proof}

\section{Framing $p$-compact flag varieties}

Before proving Theorem \ref{pittiethm}, we need an alternative description of the Pon\-trya\-gin-Thom construction \eqref{tpconstr} on $G/T$.

\begin{lemma} \label{flagmapdescription}
The map $[G/T]$ is $G$-equivariantly homotopic to the map 
\[
S_G \xrightarrow{\operatorname{incl}} BG^\g \xrightarrow{\tau} BT^\t \simeq \Sigma^r \spp{BT} \xrightarrow{\Sigma^r \epsilon} \S^r,
\]
where $BG^\g$, $BT^\t$, and $\tau$ are as in Section \ref{transfersplitting}, and all spectra except $S_G$ have a trivial $G$-action.
\end{lemma}
\begin{proof}
Applying homotopy $G$-orbits to \eqref{tpconstr}, we get a $G$-equivariant diagram
\[
\xymatrix{
S_G \ar[r] \ar[d]^{\operatorname{incl}} & S_G \sm D(\spp{G/T}) \ar[r]^-{\sim} & \spp{G/T} \sm S_T \ar[r] \ar[d] & S_T \ar@{=}[d]\\
BG^\g \ar[r]^{\tau} & BT^\t \ar[r]^-{\sim} & \spp{BT} \sm S_T \ar[r] & S_T\\
}
\]
which is commutative by the definition of $\tau$ \cite[Def. 25]{bauer:pcfm}.
\end{proof}

In the proof of Theorem \ref{pittiethm}, certain special subgroups will play an important role. In order to define them we need to recall certain facts about the Weyl group of a $p$-compact group.

Dwyer and Wilkerson showed in their ground-breaking paper \cite{dw} that given any connected $p$-compact group $G$ with maximal torus $T$, there is an associated \emph{Weyl group} $W(G)$, which is defined as the group of components of the homotopy discrete 
space of automorphisms of the fibration $BT \to BG$. This generalizes the notion of Weyl groups of compact Lie groups; they are canonically subgroups of $\GL(H_1(T)) = \GL_r(\Z_p)$, and they are so-called finite complex reflection groups. This 
means that they are generated by elements (called reflections or, more classically, pseudo-reflections) that fix hyperplanes in $\Z_p^r$. The complete classification of complex reflection groups over $\CC$ is classical and due to Shephard and Todd 
\cite{st}, the refinement to the $p$-adics is due to Clark and Ewing \cite{clark-ewing}.


Call a reflection $s \in W$ \emph{primitive} if there is no reflection $s' \in W$ of strictly larger order such that $s = (s')^k$ for some $k$.

Denote by $s \in W$ a primitive reflection of minimal order $l>1$. Let $T^s < T$ be the fixed point subtorus under $s$. Since $s$ is primitive, 
\[
\langle s \rangle = \{ w \in W \mid w|_{T^s} = \id_{T^s} \}.
\]

\begin{defn}\label{parabolicrankone}
Given connected $p$-compact group $G$ and a primitive reflection $s\in W(G)$ of minimal order $l>1$, define $C_s$ to be the centralizer of $T^s$ in $G$.
\end{defn}

Since $G$ is connected, so is the subgroup $C_s$ \cite[Lemma 7.8]{dwyer-wilkerson:center}.
Furthermore, $C_s$ has maximal rank because $T<C_s$ by definition, and the inclusion $C_s < G$ induces the inclusion of Weyl groups $\langle s \rangle < W$  \cite[Thm. 7.6]{dwyer-wilkerson:center}. Since the Weyl group of $C_s$ is $\Z/l$, the 
quotient of $C_s$ by its $p$-compact center, $C_s/Z(C_s)$, can have rank at most $1$. By the (almost trivial) classification of rank-$1$ $p$-compact groups, we find that its rank is equal to $1$ and
\begin{equation}\label{whatish}
C_s \cong \left(L_p(\S^1)^{r-1} \times L_p\S^{2l-1}\right) / \Gamma,
\end{equation}
where $L_p\S^{2l-1}$ is simply $L_p\SU(2)$ for $l=2$, and the Sullivan group given by 
\[
L_p \S^{2l-1} = \Omega L_p \left(L_p(B\S^1)_{h\Z/l}\right)
\]
for $p$ odd, and $\Gamma$ is a finite central subgroup.

\begin{proof}[Proof of Thm. \ref{pittiethm}]
By Lemma \ref{flagmapdescription}, showing equivariant null-homotopy is equivalent to showing that the map
\[
h(G/T)\colon BG^\g \xrightarrow{\tau} BT^\t \xrightarrow{\operatorname{proj}} \S^r
\]
is null. Note that for any given subgroup $H<G$ of maximal rank, there is a factorization of $\tau$ through $BH^\h$. In particular, we may assume that $G$ is connected. By the dimension hypothesis of the theorem, $W(G)$ is nontrivial. If $H=C_s$ is the subgroup associated to a primitive reflection $s\in W(G)$ of minimal order $l>1$, then the map $h(C_s/T)$ is the $(r-1)$-fold suspension of $h(L_p\S^{2l-1}/S)$ by \eqref{whatish}. Therefore, it is enough to prove the theorem for those $p$-compact groups $C_s$.

We distinguish two cases.

First suppose that $l=2$. By the classification of complex reflection groups \cite{st}, and with the terminology of that paper, this is always the case except when $W$ is a product of any number of groups from the list
\[
\{ G_4,G_5,G_{16},G_{18},G_{20},G_{25},G_{32}\}.
\]
This comment is only meant to intimidate the reader and is insubstantial for what follows.


In this case, the map $h(L_p \SU(2)/S)$ is null by Pittie \cite{Pittie:flag,Pittie-Smith:flag} since the spaces involved are Lie groups, thus $h(G/T)\simeq *$.

Now suppose that $l>2$. This forces $p>2$ as well, and since $\langle s \rangle$ acts faithfully on some line in $H^1(T;\Z_p)$ while fixing the complementary hyperplane, we must have that it acts by an $l$th root of unity, and thus $l \mid p-1$. 
The proof is finished if we can show that
\[
h(L_p\S^{2l-1}/S) = 0
\]
where $S$ is the $1$-dimensional maximal torus in the Sullivan group $G' := L_p\S^{2l-1}$.
To see this, note that the inclusion $S \to G'$ factors through the maximal torus normalizer $N_{G'}(S) \cong S \rtimes \Z/l$, and thus
\[
h(G'/S)\colon B{G'}^{\g'} \xrightarrow{\tau_1} BN^\n \xrightarrow{\tau_2} \Sigma\spp{BS} \to \S^1.
\]
If $P \simeq \bigvee_{i=0}^{p-2} e_i P$ is any stable splitting of the $p$-completed complex projective space $P=\spp{BS}$ induced by idempotents $e_i$ as in the previous section, then the rightmost projection map clearly factors through $e_0P$, 
which is the part containing the bottom cell. Since $p>2$, Corollary \ref{transferfactorization} shows that there is an idempotent $f\in [P,P]$ such that $\tau \simeq f \circ \tau$ and $\tau \circ e_0 = e_0 \circ \tau = 0$, proving the theorem.
\end{proof}

\section{The adjoint representation}

Let $G$ be a $d$-dimensional connected $p$-compact group with maximal torus $T$ of rank $r$. Choose a set $\{s_1,\dots,s_{r'}\}$ of generating reflections of $W=W(G)$ with $r'$ minimal. The classification of pseudo-reflection groups \cite{st,clark-ewing} 
implies that for $G$ semisimple, most of the time $r=r'$, but there are cases where $r'=r+1$.

\begin{example}[The group no. 7]
Let $p \equiv 1 \pmod{12}$. Let $G_7$ be the finite group generated by the reflection $s$ of order $2$ and the two reflections $t,\;u$ of order $3$, where $s,\;t,\;u \in GL_2(\Z_p)$ are given by
\[
s = \begin{pmatrix} -1 & 0 \\ 0 & 1 \end{pmatrix}, \quad
t = \frac1{\sqrt{2}} \begin{pmatrix} -\zeta & \zeta^7 \\ -\zeta & -\zeta^7 \end{pmatrix},\quad
u = \frac1{\sqrt{2}} \begin{pmatrix} -\zeta^7 & -\zeta^7 \\ \zeta & -\zeta \end{pmatrix}.
\]
Here $\zeta$ is a $24$th primitive root of unity. Note that although possibly $\zeta \not\in \Z_p$, $\frac1{\sqrt{2}} \zeta \in \Z_p$.
In Shephard and Todd's classification, this is the restriction to $\Z_p$ of the complex pseudo-reflection group no. 7. They show that even over $\CC$, $G_7$ cannot be generated by two reflections. The associated $p$-compact group is given by
\[
\Omega L_p((BT^2)_{hG_7}).
\]
\end{example}

If $G$ is not semisimple (i.~e. it contains a nontrivial normal torus subgroup), then $r'$ may be smaller than $r$. Set $\kappa = r+1-r' \geq 0$.

Let $\I_{r'}$ be the set of proper subsets of $\{1,\dots,r'\}$, and for $I \subseteq \{1,\dots,r'\}$, let $T_I$ be the fixed point subtorus $T^{\langle s_i \mid i \in I \rangle}$ and $C_I = C_G(T_I)$ be the centralizer in $G$, which is connected by \cite[Lemma 7.8]{dwyer-wilkerson:center}.

\begin{defn}
Let $G$ be a connected $p$-compact group. Define the \emph{adjoint space} $A_G$ by the homotopy colimit
\[
A_G = \Sigma^\kappa \hocolim_{I \in \I_r} G/C_I
\]
with the induced left $G$-action, and the trivial $G$-action on the suspension coordinates.
\end{defn}

Theorem \ref{ckthm} shows that if $G$ is a the $p$-completion of a connected, semisimple Lie group (in this case $r=r'$ and $\kappa=1$), then $A_G$ is a $d$-dimensional sphere $G$-equivariantly homotopy equivalent to $\g \cup \{\infty\}$. This holds 
more generally: if $G$ is a connected, compact Lie group with maximal normal torus $T^k$ then
\[
A_G \cong \Sigma^k A_{G/T^k} = (\t\cup \{\infty\}) \sm (\g/\t \cup \{\infty\}) = \g \cup \{\infty\}.
\]

\begin{lemma} \label{mayervietorislemma}
Let $\I_{r}$ be the poset category of proper subsets of $\{1,\dots,r\}$ and 
\[
F,G \colon \I_k \to \{\text{finite CW-complexes or finite CW-spectra}\}.
\]
be two functors. Then
\begin{enumerate}
\item If $F$ is has the property that $\dim F(\emptyset) > \dim F(I)$ for every $I \neq \emptyset$, then
\[
\dim \hocolim F = \dim F(\emptyset)+k-1.
\]
\item If $f\colon F \to G$ is a natural transformation of two such functors such that
\[
f_*(\emptyset) \colon H_{\dim F(\emptyset)}(F(\emptyset)) \xrightarrow{\cong} H_{\dim G(\emptyset)}(G(\emptyset)),
\]
then f induces an isomorphism
\[
\hocolim f_*\colon H_{\dim \hocolim F}(\hocolim F) \to H_{\dim \hocolim G}(\hocolim G).
\]
\item Let $F\colon \I_r \to \Top$ be the functor given by $F(\emptyset) = \S^n$, $F(I)=*$ for $I \neq \emptyset$. Then $\hocolim_{\I_r} F \simeq \S^{n+r-1}$.
\end{enumerate}
\end{lemma}
\begin{proof}
The first two assertions follow from the Mayer-Vietoris spectral sequence \cite[Chapter XII.5]{yellowmonster},
\[
E^1_{p,q} = \smashop\bigoplus_{I \in \I_k,\;|I|=k-1-p} H_q(F(I)) \Longrightarrow H_{p+q}(\hocolim F),
\]
along with the observation that under the dimension assumptions of (1), $E^1_{p,q}=0$ for $q \geq \dim F(\emptyset)$ except for $E^1_{k-1,\dim F(\emptyset)} = H_{\dim F(\emptyset)}(F(\emptyset))$. In particular, this group cannot support a nonzero 
differential and thus
\[
H_i(F(\emptyset)) \cong H_{i+k-1}(\hocolim F) \quad \text{for } i \geq \dim F(\emptyset).
\]
The third one is an immediate consequence of the Mayer-Vietoris spectral sequence.
\end{proof}

\begin{corollary}
For any connected $p$-compact group $G$, $A_G$ is a $d$-dimensional $G$-space.
\end{corollary}
\begin{proof}
This follows from Lemma \ref{mayervietorislemma}. Indeed, since any $C_I$ ($I \neq \emptyset$) is connected and has the nontrivial Weyl group $W_I$, its dimension is greater than $\dim T$. So the condition
\[
\dim F(\emptyset) = \dim G/T > \dim F(I)
\]
is satisfied, and
\[
\dim \hocolim F = d-r+r'-1 = d-\kappa.
\]
\end{proof}

As mentioned at the end of the introduction, for $p$-compact groups $G$, $A_G$ is not usually a sphere, as the following example illustrates.

\begin{example}
Let $p \geq 5$ be a prime, and let $G = \S^{2p-3}$ be the Sullivan sphere, whose group structure is given by $BG = L_p \left(BS_{hC_{p-1}}\right)$, where $C_{p-1} \subseteq \Z_p^\times$ acts on $BS = K(\Z_p,2)$ by multiplication on $\Z_p$. Clearly, 
$G$ has rank $1$, and $\I_1$ consists only of a point, thus $A_G = \Sigma G/T \simeq L_p \Sigma \CP^{p-2}$. Since $p \geq 5$, this is not a sphere.
\end{example}

For the proof of Theorem~\ref{adjointcomparison} we need a preparatory result.

\begin{prop} \label{simplenullhomotopy}
Let $P$ be a $p$-compact subgroup of maximal rank in a $p$-compact group $G$. Denote by $T$ a maximal torus of $P$ (and thus also of $G$). Then the following composition is $G$-equivariantly null-homotopic:
\[
f_{G,P}\colon S_G \sm DS_T \to \spp{G/T} \to \spp{G/P}.
\]
The second map is the canonical projection, whereas the first map is given by using the duality isomorphism
\begin{multline*}
S_G \sm DS_T \to S_G \sm D(\spp{G/T}) \sm DS_T \xrightarrow{\sim} G_+ \sm_T S_T \sm DS_T \\
\simeq \spp{G/T} \sm S_T \sm DS_T \xrightarrow{\id \sm \operatorname{ev}} \spp{G/T}.
\end{multline*}
\end{prop}

\begin{proof}
In \cite[Cor. 24]{bauer:pcfm} it was shown that the relative duality isomorphism from Lemma \ref{relativeduality} is natural in the sense that the following diagram commutes:
\[
\xymatrix{
(\spp{G})^{hP^{\op}} \ar[d]^{\operatorname{res}} \ar[r]^{\sim} & G_+ \sm_P S_P \ar[r]^{\sim} & D(\spp{G/P}) \sm S_G \ar[d]^{D(\text{proj}) \sm \id} \\
(\spp{G})^{hT^{\op}}\ar[r]^{\sim} & G_+ \sm_T S_T \ar[r]^{\sim} & D(\spp{G/T}) \sm S_G
}
\]
Taking duals and smashing with $DS_G$, we find that the map of the proposition is the left hand column in the diagram
\[
\xymatrix{
S_G \sm D(G_+ \sm_P S_P) \ar[r]^-{\sim} & \spp{G/P}\\
S_G \sm D(G_+ \sm_T S_T) \ar[u] \ar[r]^-{\sim} & \spp{G/T} \ar[u]\\
S_G \sm D(\spp{G/T}) \sm DS_T \ar[u]^{\sim}\\
S_G \sm DS_T \ar[u]
}
\]
Thus we need to show that the composition
\[
G_+ \sm_P S_P \to G_+ \sm_T S_T \simeq \spp{G/T} \sm S_T \to S_T
\]
is $G$-equivariantly trivial, or equivalently, that
\[
S_P \to P_+ \sm_T S_T \to S_T
\]
is $P$-equivariantly trivial. But this map is exactly the homotopy class represented by $[P/T]$, thus the assertion follows from Theorem \ref{pittiethm}.
\end{proof}


\begin{proof}[Proof of Thm. \ref{adjointcomparison}]
Let $G$ be a connected $p$-compact group whose Weyl group is generated by a minimal set of $r'$ reflections.
Let $F,A\colon \I_{r'} \to \Top$ be the functors given by $F(\emptyset) = S_G \sm DS_T$, $F(I)=*$ for $I \neq \emptyset$, and $A(I)=G/C_I$. Note that, since $G$ is connected, $C_G(T)=T$ \cite[Proposition 9.1]{dw} and $A(\emptyset)=G/T$. There is a map $\Phi\colon F \to A$ of $\I_{r'}$-diagrams in the homotopy category of $G$-spectra which is fully described by defining 
\[
\Phi(\emptyset)=f_{G,T}\colon F(\emptyset) = S_G \sm DS_T \to \spp{G/T}
\]
as the map given in Prop.~\ref{simplenullhomotopy}.
The strategy of the proof is to obtain a functor $F \colon \I_{r'} \to \Top$ such that $F(\emptyset) = S_G \sm DS_T$, $F(I)\simeq *$ for $I \neq \emptyset$, and  a map of $\I_{r'}$-diagrams $\Phi \colon F \to A$ in the category of $G$-spectra such that $\Phi(\emptyset) = f_{G,T}$. From this we get a $G$-equivariant map
\[
S_G \simeq \S^\kappa \sm \Sigma^{r'-1} S_G \sm DS_T \simeq \S^\kappa \sm \hocolim_{\I_{r'}} F \to \Sigma^\kappa \hocolim_{\I_{r'}} \spp{G/C_I} \simeq \spp{A_G},
\]
which will give us the splitting.


We will proceed by induction on the number of generating reflections $r'$. If $r'=1$ then $A_G$ is $\S^\kappa \sm G/T$ and $\Phi(\emptyset)=\S^\kappa \sm f_{G,T}$.  
 We can construct the functor $F$ and the natural transformation $\Phi$ step by step. Fix a subset $I$ of cardinality $k$, and assume that $F$ and $\Phi$ have been defined for all vertices in the diagram corresponding to $I'$ with $|I'|<k$. 

Let $\mathcal P(I)$ be the poset category of all proper subsets of $I$. Since $F$ and $\Phi$ are defined over $\mathcal P(I)$ by induction hypothesis, we can consider $\hocolim_{\mathcal P(I)} F\simeq \Sigma^{k-1} S_G \sm DS_T \to \spp{G/C_I}$. It is enough to show that this map is 
$G$-equivariantly nullhomotopic. Then, we can fix a null-homotopy and extend the map to the cone of $\hocolim_{\mathcal P(I)} F$. Finally, we define $F(I)=C(\hocolim_{\mathcal P(I)} F)$ and $\Phi(I)$ is the corresponding extension.

Note that $S_G \sm DS_T \to \spp{G/T}$ factors through $G_+ \sm_{C_I} S_{C_I} \sm DS_T$. By induction, we know there is a map 
\[
 \Sigma^{k-1} S_{C_I} \sm DS_T \to  \spp{\hocolim_{J \in \I_{k}}{C_I}/C_J},
\]
which splits the top cell. We get a factorization
\begin{multline}
 \Sigma^{k-1} S_{G} \sm DS_T \to \Sigma^{k-1} G_+\sm_{C_I} S_{C_I} \sm DS_T \\ \to  G_+\sm_{C_I} \spp{\hocolim_{J \in \I_{k}}{C_I}/C_J} \to G_+\sm_{C_I} \S^0.
\end{multline}
It thus suffices to show that in the $\I_k$-diagram
\[
\xymatrix{
S_{C_I} \sm DS_T \ar[d] \ar[r] \ar@<1mm>[r] \ar@<-1mm>[r] & {\{ * \}_{J \in \I_{k}-\{\emptyset\}}} \ar[r]^-{\hocolim} \ar[d] & \Sigma^{k-1} S_{C_I} \sm DS_T \ar[d]\\
{\spp{{C_I}/T}} \ar[d] \ar[r] \ar@<1mm>[r] \ar@<-1mm>[r]  &  {\{ \spp{{C_I}/C_J} \}_{J \in \I_{k}-\{\emptyset\}}} \ar[d] \ar[r]^-{\hocolim} & \spp{\hocolim_{J \in \I_{k}}{C_I}/C_J} \ar[d]\\
{\S^0} \ar[r] \ar@<1mm>[r] \ar@<-1mm>[r]  & {\{ \S^0\}_{I \in \I_{k}-\{\emptyset\}}} \ar[r]^-{\hocolim} & \S^0
}
\]
the right hand side composition $\Sigma^{k-1} S_{C_I} \sm DS_T \to \S^0$ is ${C_I}$-equivariantly null-homotopic. In the latter diagram, it makes no difference whether the centralizers are taken in ${C_I}$ or in $G$. But by Theorem \ref{pittiethm}, 
the left hand column is already null-homotopic, thus, as a colimit of null-homotopic maps over a contractible diagram, so is the right hand column.
\end{proof}

\subsection*{Conclusion and questions}

In this paper, we have compared two imperfect notions of adjoint representations of a $p$-compact group $G$. One ($S_G$) is a sphere, but has a $G$-action only stably; the other ($A_G$) is an unstable $G$-space, but fails to be a sphere. The 
question remains whether there is an unstable $G$-sphere whose suspension spectrum is $S_G$. It might even be true that $A_G$ splits off its top cell after only one suspension, yielding a solution to this problem in the cases where the Weyl group 
of the rank-$r$ group $G$ is generated by $r$ reflections.

There are also a number of interesting open questions about the flag variety $G/T$ of a $p$-compact groups:
\begin{itemize}
\item By the classification of $p$-compact groups, $H^*(G/T;\Z_p)$ is torsion free and generated in degree $2$. Can this be seen directly?
\item Is there a manifold $M$ such that $L_pM \simeq G/T$, analogous to smoothings of $G$ \cite{BKNP:02,bauer-pedersen:localloop}? Is it a boundary of a manifold?
\item If such a manifold $M$ exists, can it be given a complex structure?
\end{itemize}

\end{document}